\documentclass[a4paper,12pt]{extarticle}

\usepackage[T1]{fontenc}	% Saisie en
\usepackage[utf8]{inputenc}	% 	franÃ§ais

\usepackage{geometry}		% Taille de la feuille
\usepackage{fancyhdr}		% Titre courant
\usepackage{setspace}		% Interligne
\usepackage{lscape}		% Mode paysage
\usepackage{multicol}		% Plusieurs colonnes
\usepackage{varioref}	% RÃ©fÃ©rences
\usepackage{float, caption}

\usepackage{calc}		% Outils de calcul
\usepackage{ifthen}		% Tests if/then/else
\usepackage{xspace}		% Gestion des espaces

\usepackage{pifont}		% Police \ding
\usepackage{eurosym}		% Symbole de l'euro
\usepackage{soul}		% Souligner
\usepackage{enumerate}		% Listes
\usepackage{verbatim}		% Codes source
\usepackage{moreverb}
\usepackage{fancybox}
\usepackage{boxedminipage}		%	et listings

\usepackage{array}		% Outils supplÃ©mentaires
\usepackage{multirow}		% Colonnes multiples
\usepackage{tabularx}		% Largeur totale donnÃ©e
\usepackage{longtable}		% Sur plusieurs pages

\usepackage{graphicx}		% Insertion de figures
\usepackage{epic}		% CapacitÃ©s graphiques
\usepackage{eepic}		% 	Ã©tendues
\usepackage{afterpage}		% Voir page 69
\usepackage{rotating}		% Tourner du texte
\usepackage{pgf,tikz}
\usepackage{mathrsfs}
\usetikzlibrary{arrows}

\usepackage{amsmath}		% Commandes essentielles
\usepackage{amssymb}		% Principaux symboles
\usepackage{mathrsfs}		% Police calligraphique
\usepackage[thmmarks,amsmath]{ntheorem}		% ThÃ©orÃ¨mes
\usepackage{dsfont}
\usepackage{youngtab}

\usepackage{pgf,tikz}
\usepackage{mathrsfs}
\usetikzlibrary{arrows}

\newtheorem{prop}{Proposition}[subsection]

\newtheorem{theo}[prop]{Theorem}
\newtheorem{cor}[prop]{Corollary}
\newtheorem{lem}[prop]{Lemma}
\newtheorem{ex}[prop]{Example}

\newtheorem{definition}[prop]{Definition}

\title{A common limit in large rank for Markov chains defined from representations of classical Lie algebras}
\author{Vivien Despax\thanks{Laboratoire de Math\'ematiques et Physique Th\'eorique, UMR-CNRS 7350,
F\'ed\'eration Denis Poisson, FR-CNRS 2964,
Universit\'e Fran\c{c}ois--Rabelais de Tours, Parc de Grandmont, 37200 Tours, France, vivien.despax@lmpt.univ-tours.fr}}

\date{}

\begin{document}
\maketitle
\begin{abstract} 
From the datum of an integer partition and a classical Lie algebra, one can define a Markov chain on an associated multiplicative graph. For each classical family $A,C,B,D$, we thus obtain a sequence of Markov chain which is indexed by the rank of the considered algebra.
%We show that each such sequence admits a limit when the rank tends to infinity. Moreover, we establish that the four limit Markov chains have the same law.
In this article we show that, for each type, the transition kernel of the Markov chain has a limit when the rank tends to infinity. Moreover, the limit kernel does not depend on the considered type. 

\end{abstract}

\section{Introduction}
In what follows, $\delta$ is a nonzero partition. For all $r$ sufficiently large, $\delta$ is regarded as a dominant weight of $\mathfrak{g}_r=\mathfrak{gl}_r, \mathfrak{sp}_{2r} ,\mathfrak{so}_{2r+1} $ or $\mathfrak{so}_{2r}$. For fixed $r$, it is then possible to define a graph $\mathcal{G}^{\left(\delta,\mathfrak{g}_r\right)}$ reflecting the multiplicative structure that arises when one decomposes the successive tensor powers \[ V^{\mathfrak{g}_r}\left( \delta\right),V^{\mathfrak{g}_r}\left( \delta\right)^{\otimes 2},V^{\mathfrak{g}_r}\left( \delta\right)^{\otimes 3},\ldots\]into irreducible components. Introducing then a stochastic matrix $\Pi^{\left(\delta,\mathfrak{g}_r\right)}$ on $\mathcal{G}^{\left( \delta,\mathfrak{g}_r\right)}$, we now have a Markov chain $\mathcal{M}^{\left(\delta,\mathfrak{g}_r \right)}$. The aim of our article is to study what happens when $r$ tends to infinity.

The Markov chain $\mathcal{M}^{\left(\delta,\mathfrak{g}_r \right)}$ originally appears in~\cite{O,LLP1,LLP2} as a random walk on the weight lattice of $\mathfrak{g}_r$ conditioned to never exit the cone of dominant weights. In these articles, it is shown that the transition matrix of such conditioned process can be expressed in terms of specializations of the corresponding normalized Weyl characters. The determination of the harmonic functions on the graph $\mathcal{G}^{\left( \delta,\mathfrak{g}_r\right)}$ is a closely connected subject which is studied in~\cite{LT}. 
It is then natural to ask what happens when $r$ tends to infinity.

Roughly speaking, we prove in this article that the four classical families \[A=\left(\mathfrak{gl}_r \right)_r,C= \left(\mathfrak{sp}_{2r} \right)_r, B=\left(\mathfrak{so}_{2r+1} \right)_r, D=\left(\mathfrak{so}_{2r} \right)_r\] give rise to the same limit Markov chain. Our result intersects with a similar phenomenon observed by Nikitin and Vershik in~\cite{NV} for some graded graphs and their Pascalized versions. Nevertheless, the graphs that we consider in this article cannot be seen in general as Pascalized versions of simpler graded graphs. 

To give sense to this statement, we start by observing that the sequence of multiplicative graphs $\left(\mathcal{G}^{\left(\delta,\mathfrak{g}_r\right)}\right)_r$ admits a natural limit. This is a consequence of the stabilization of the tensor multiplicities for large rank. Next, we prove that the sequence of stochastic matrices $\left(\Pi^{\left(\delta,\mathfrak{g}_r\right)}\right)_r$ is convergent. To accomplish this, we establish the result in type $A$ by showing that the specializations of the normalized Schur functions used to define $\Pi^{\left(\delta,\mathfrak{gl}_r\right)}$ converge when $r$ tends to
infinity. This mainly follows from the Weyl character formula. By using the notion of Kashiwara-Nakashima tableaux, we show the convergence to the same limit in type $C,B,D$. 
This stabilization in large rank and the existence of a common limit when this
rank tends to infinity illustrates a general phenomenon in the combinatorics
of classical root systems. For example each tensor product multiplicity
$c_{\delta,\kappa}^{\lambda}(\mathfrak{g}_{r})$ appearing in the decomposition%
\[
V^{\mathfrak{g}_{r}}(\delta)\otimes V^{\mathfrak{g}_{r}}(\kappa)=%
%TCIMACRO{\tbigoplus \limits_{\lambda}}%
%BeginExpansion
{\textstyle\bigoplus\limits_{\lambda}}
%EndExpansion
V^{\mathfrak{g}_{r}}(\lambda)^{c_{\delta,\kappa}^{\lambda}(\mathfrak{g}_{r})}%
\]
stabilizes in large rank and the common limit when $r$ tends to infinity is
the Littlewood-Richardson coefficients $c_{\delta,\kappa}^{\lambda}$ (see for
example \cite{L}).

The article is organized as follows. Section~\ref{b} introduces the required material and notation on representation theory. %(§~\ref{rtco})
Section~\ref{fr} describes the Markov chain $\mathcal{M}^{\left( \delta,\mathfrak{g}_r\right)}$ on finite rank multiplicative graph. %(§~\ref{fr}).
Section~\ref{irmc} is concerned with the limit Markov chains. We first define the limit graph (§~\ref{irg}) and then study the asymptotic behaviour of $\left( \Pi^{\left(\delta,\mathfrak{g}_r\right)}\right)_r$ (§~\ref{irm}). 
We start with the type $A$ case (§~\ref{A}) and pursue with types $C,B,D$ (§~\ref{CBD}).
Our main theorem is then stated in §~\ref{main}. Finally in Section \ref{SecFR} we link our results to the generalized Pitman transform defined in \cite{BBOC1} and to the study of the extremal harmonic functions on multiplicative graphs achieved in \cite{LT}.

\section{Background on representation theory and combinatorics}
\label{b}
In this article, we only consider complex Lie algebras and their finite-dimensional representations. Let $X=\left(\mathfrak{g}_r\right)_r$ be one of the four classical families $A=\left(\mathfrak{gl}_r \right)_r,C= \left(\mathfrak{sp}_{2r} \right)_r, B=\left(\mathfrak{so}_{2r+1} \right)_r$ or $D=\left(\mathfrak{so}_{2r} \right)_r$.
% \subsection{Representation theory and combinatorics}
% \label{rtco}
Let us fix a positive integer $r$ for the whole section. We introduce some objects defined from $\mathfrak{g}_r$. Nevertheless, we do not include this parameter in our notation for the moment.
 
\subsection{Root system and weight lattice}
The Cartan subalgebra of $\mathfrak{g}_r$ consisting of the diagonal matrices is denoted by $\mathfrak{h}$. There exists an Euclidean space $\left(E,\left<\,.\, ,\,.\,\right> \right)$ endowed with an orthonormal basis $\left(\varepsilon_1,\varepsilon_2,\ldots,\varepsilon_r \right)$ such that the root system $R$ of $\mathfrak{g}_r$ (with respect to $\mathfrak{h}$) is identified to the following finite subset of $E$:\[R=\begin{cases}\left\lbrace \pm\left( \varepsilon_i-\varepsilon_j \right) : 1\leq i<j\leq r\right\rbrace &\text{if }\mathfrak{g}_r=\mathfrak{gl}_r\\
\left\lbrace \pm 2\varepsilon_i : 1\leq i \leq r \right\rbrace \cup \left\lbrace  \pm\left( \varepsilon_i\pm\varepsilon_j \right)  : 1\leq i<j\leq r\right\rbrace &\text{if }\mathfrak{g}_r=\mathfrak{sp}_{2r}\\
\left\lbrace \pm \varepsilon_i : 1\leq i \leq r \right\rbrace \cup \left\lbrace \pm\left( \varepsilon_i\pm\varepsilon_j\right)  : 1\leq i<j\leq r\right\rbrace &\text{if }\mathfrak{g}_r=\mathfrak{so}_{2r+1}\\
\left\lbrace  \pm\left( \varepsilon_i\pm\varepsilon_j\right)   : 1\leq i<j\leq r\right\rbrace &\text{if }\mathfrak{g}_r=\mathfrak{so}_{2r}\\
\end{cases}.\]
Denote by $F$ the subspace of $E$ spanned by $R$. The Weyl group $W$ of $R$ is regarded here as the subgroup of $\mathrm{O}\left( F\right) $ generated by the orthogonal reflections with fixed points the $\left\lbrace \alpha\right\rbrace^{\perp} ,\alpha \in R$. Let $I$ be $\left\lbrace 1,2,\ldots,r-1\right\rbrace $ if $\mathfrak{g}_r=\mathfrak{gl}_{r}$ and $\left\lbrace 1,2,\ldots,r\right\rbrace $ in the other cases. We introduce the simple system $S=\left(\alpha_i \right) _{i\in I}$ such that $\alpha_i=\varepsilon_i-\varepsilon_{i+1}$ if $1\leq i<r$ and \[\alpha_r=\begin{cases} 2\varepsilon_r   &\text{if }\mathfrak{g}_r=\mathfrak{sp}_{2r}\\
\varepsilon_r   &\text{if }\mathfrak{g}_r=\mathfrak{so}_{2r+1}\\
\varepsilon_{r-1}+\varepsilon_r    &\text{if }\mathfrak{g}_r=\mathfrak{so}_{2r}\\
\end{cases}\]
which yields the partition $R=R_+ \sqcup \left( -R_+\right) $ where\[R_+=\begin{cases} \left\lbrace \varepsilon_i-\varepsilon_j : 1\leq i<j\leq r\right\rbrace &\text{if }\mathfrak{g}_r=\mathfrak{gl}_r\\
\left\lbrace 2\varepsilon_i : 1\leq i \leq r \right\rbrace \cup \left\lbrace  \varepsilon_i\pm\varepsilon_j  : 1\leq i<j\leq r\right\rbrace &\text{if }\mathfrak{g}_r=\mathfrak{sp}_{2r}\\
\left\lbrace  \varepsilon_i : 1\leq i \leq r \right\rbrace \cup \left\lbrace \varepsilon_i\pm\varepsilon_j : 1\leq i<j\leq r\right\rbrace &\text{if }\mathfrak{g}_r=\mathfrak{so}_{2r+1}\\
\left\lbrace  \varepsilon_i\pm\varepsilon_j  : 1\leq i<j\leq r\right\rbrace &\text{if }\mathfrak{g}_r=\mathfrak{so}_{2r}\\\end{cases}.\]
The family of fundamental weights $\Omega=\left(\omega_i \right)_{i\in I}$ associated to $S$ is such that
\[\begin{cases} \omega_i=\sum_{j=1}^i \varepsilon_j,i\in I&\text{if }\mathfrak{g}_r=\mathfrak{gl}_r,\mathfrak{sp}_{2r}\\
\omega_i=\sum_{j=1}^i \varepsilon_j, 1\leq i<r,\omega_r=\frac{1}{2}\sum_{j=1}^r \varepsilon_j&\text{if }\mathfrak{g}_r=\mathfrak{so}_{2r+1}\\
\omega_i =\sum_{j=1}^i \varepsilon_j,1\leq i<r-1, \omega_{r-1}=\frac{1}{2}\left(\sum_{j=1}^{r-1} \varepsilon_j\right) -\frac{1}{2}\varepsilon_r,\omega_r=\frac{1}{2}\sum_{j=1}^r \varepsilon_j &\text{if }\mathfrak{g}_r=\mathfrak{so}_{2r}
\end{cases}.\]
If we set $\omega_i^{\vee}=\frac{2}{\left< \alpha_i,\alpha_i\right> }\omega_i$, then the family $\Omega^{\vee}=\left( \omega_i^{\vee}\right) _{i\in I}$ is such that\[\left< \omega_i^{\vee},\alpha_j\right> =\delta_{i\,j}\qquad i,j\in I.\]
Notice that $\omega_i^{\vee}=\omega_i$ when $X=A$.
Introduce the weight lattice of $\mathfrak{g}_r$
\[L=\begin{cases}\bigoplus_{i=1}^r\mathbb{Z}\varepsilon_i&\text{if }\mathfrak{g}_r=\mathfrak{gl}_r\\
\bigoplus_{i=1}^r \mathbb{Z}\omega_i=\bigoplus_{i=1}^r \mathbb{Z}\varepsilon_i&\text{if }\mathfrak{g}_r=\mathfrak{sp}_{2r}\\
\bigoplus_{i=1}^r \mathbb{Z}\omega_i=\bigoplus_{i=1}^r \mathbb{Z}\varepsilon_i+\mathbb{Z}\left( \frac{1}{2}\sum_{i=1}^{r}\varepsilon_i\right) &\text{if }\mathfrak{g}_r=\mathfrak{so}_{2r+1},\mathfrak{so}_{2r}\\
\end{cases}\]
and set\[P=\begin{cases}\bigoplus_{i=1}^r \mathbb{Z}_+\varepsilon_i&\text{if }\mathfrak{g}_r=\mathfrak{gl}_r\\
L&\text{if }\mathfrak{g}_r=\mathfrak{sp}_{2r},\mathfrak{so}_{2r+1},\mathfrak{so}_{2r}.
\end{cases}\]
The weights of the (polynomial if $\mathfrak{g}_r=\mathfrak{gl}_r$) representations of $\mathfrak{g}_r$ are identified with the elements of $P$: if $V$ is a representation of $\mathfrak{g}_r$, then it admits a weight space decomposition\[V=\bigoplus_{\omega\in P} V_\omega=\bigoplus_{\omega\in P\left( V\right) } V_\omega\]where\[V_\omega =\left\lbrace v\in V : hv=\omega\left( h\right) v \textit{ for any $h$ in $\mathfrak{h}$}\right\rbrace\]and \[P\left( V\right)=\left\lbrace \omega\in P: V_\omega \neq \left\lbrace 0\right\rbrace \right\rbrace.\]

\subsection{Irreducible representations and characters}
Introduce the set of dominant weights of $\mathfrak{g}_r$\[P_+=\begin{cases}\left\lbrace \sum_{i=1}^r \lambda_i \varepsilon_i \in P : \lambda_1\geq \lambda_2\geq\ldots\geq \lambda_r \right\rbrace&\text{if }\mathfrak{g}_r=\mathfrak{gl}_r\\
\bigoplus_{i=1}^r \mathbb{Z}_+ \omega_i=\left\lbrace \sum_{i=1}^r \lambda_i \varepsilon_i \in P : \lambda_1\geq \lambda_2\geq\ldots\geq \lambda_r\geq 0 \right\rbrace&\text{if }\mathfrak{g}_r=\mathfrak{sp}_{2r}\\
\bigoplus_{i=1}^r \mathbb{Z}_+ \omega_i=\left\lbrace \sum_{i=1}^r \lambda_i \varepsilon_i \in P : \lambda_1\geq \lambda_2\geq\ldots\geq \lambda_r\geq 0 \right\rbrace&\text{if }\mathfrak{g}_r=\mathfrak{so}_{2r+1}\\
\bigoplus_{i=1}^r \mathbb{Z}_+ \omega_i=\left\lbrace \sum_{i=1}^r \lambda_i \varepsilon_i \in P : \lambda_1\geq \lambda_2\geq\ldots\geq \lambda_{r-1}\geq \left| \lambda_r\right|  \right\rbrace&\text{if }\mathfrak{g}_r=\mathfrak{so}_{2r}\\
\end{cases}.\]
The set $P_+$ parametrizes the irreducible (polynomial if $\mathfrak{g}_r=\mathfrak{gl}_r$) representations of $\mathfrak{g}_r$ and every representation $V$ of $\mathfrak{g}_r$ decomposes into irreducible components appearing with multiplicities
\[V=\bigoplus_{\lambda\in P_+} V\left( \lambda\right) ^{\oplus m_{V\,\lambda}}\]
where $V\left( \lambda\right) $ is the irreducible representation of ${\mathfrak{g}_r}$ corresponding to $\lambda$. In particular, for any $\delta$ in $P_+$, the tensor powers $V\left( \delta\right) ^{\otimes n}$ ($n>0$) and the tensor products $V\left( \lambda\right) \otimes V\left( \delta\right)$ ($\lambda\in P_+$) admit such decomposition:\[V\left( \delta\right) ^{\otimes n}=\bigoplus_{\mu\in P_+} V\left(\mu \right)^{\oplus f_{n\,\mu}}\qquad  V\left( \lambda\right) \otimes V\left( \delta\right) =\bigoplus_{\mu\in P_+} V\left(\mu \right)^{\oplus m_{\lambda\,\mu}}. \]

Given a list of variables $\left(x_i\right)_{1\leq i\leq r} $ and $\omega$ in $L$, we set $x^\omega =\prod_{i=1}^r x_i^{\left<\omega,\varepsilon_i \right> }$. For $\lambda$ in $P_+$, the Weyl character of the irreducible representation $V\left( \lambda\right) $ is \[s_\lambda \left( x\right) =\sum_{\omega\in P}K_{\lambda\,\omega}x^\omega=\sum_{\omega\in P\left( \lambda\right) }K_{\lambda\,\omega}x^\omega\]
where $K_{\lambda\,\omega}=\dim V\left(\lambda\right) _\omega$ and $P\left( \lambda\right)=P\left( V\left( \lambda\right) \right)$. Recall the Weyl character formula: \[s_\lambda\left(x\right)=\frac{a_{\lambda+\rho} \left({x}\right) }{a_{\rho} \left({x}\right) }\]
where\[a_\omega \left({x}\right)=\sum_{\sigma\in W} \epsilon\left( \sigma\right) {x}^{\sigma\omega}\qquad \omega\in L\]$\epsilon$ being the sign character of $W$ and $\rho=\sum_{i\in I}\omega_i$. The dominator $a_{\rho}\left( x\right) $ can be expressed as\[a_{\rho} \left( {x}\right)={x}^{\rho}\prod_{\alpha\in R_+}\left( 1-{x}^{-\alpha}\right).\]

\subsection{Partitions}
\begin{definition}
	Any weakly decreasing sequence of nonnegative integers with a finite number of nonzero terms is called a partition. 
\end{definition}
Each partition is identified to its Young diagram. Given a partition $\lambda$, we denote by $\ell\left( \lambda\right)$ the number of its parts, that is the number of its nonzero terms. We denote by $\left|\lambda\right|$ the sum of its terms, that is its number of boxes. The set of partitions with at most $\ell$ parts is denoted by $\mathcal{P}_\ell$. The sequence $\left(\mathcal{P}_\ell \right) _{\ell\in \mathbb{N}}$ is increasing and the set of all partitions is $\mathcal{P}_\infty=\bigcup_{\ell\in \mathbb{N}} \mathcal{P}_\ell$. If $\lambda=\left(\lambda_1,\lambda_2,\ldots,\lambda_\ell,0,0,\ldots \right)$ is in $\mathcal{P}_\ell$ with $\ell\leq r$, then we identify $\lambda$ with the element  $\lambda_1\varepsilon_1+\lambda_2\varepsilon_2+\ldots+\lambda_\ell\varepsilon_\ell$ of $P_+$. According to this, one has then $\mathcal{P}_r=P_+$ if $\mathfrak{g}_r=\mathfrak{gl}_r,\mathfrak{sp}_{2r}$ and $\mathcal{P}_r\subset P_+$ if $\mathfrak{g}_r=\mathfrak{so}_{2r+1},\mathfrak{so}_{2r}$.

If $\delta$ and $\lambda$ are in $\mathcal{P}_r$, it can be shown that we have in fact\begin{equation}V\left( \delta\right) ^{\otimes n}=\bigoplus_{\mu\in \mathcal{P}_r} V\left(\mu \right)^{\oplus f_{n\,\mu}}\qquad  V\left( \lambda\right) \otimes V\left( \delta\right) =\bigoplus_{\mu\in \mathcal{P}_r} V\left(\mu \right)^{\oplus m_{\lambda\,\mu}}.\label{pie}\end{equation} We can be even more precise by showing that we can restrict the above direct sums to $\mu$ in $\mathcal{P}_r$ such that $\left| \mu\right|\leq n \left|\delta \right|$ and $\left| \mu\right|\leq \left|\lambda \right| + \left|\delta \right|$, respectively.

Given $\lambda$ in $\mathcal{P}_r$, the representation $V\left(\lambda\right)$ appears with a positive multiplicity in the decomposition of $V\left(\tiny{\yng(1)}\right)^{\otimes \left| \lambda\right|}$ into irreducible. One deduces in particular\[\dim V\left(\lambda\right)\leq \left(\dim V\left(\tiny{\yng(1)}\right)\right)^{\left| \lambda\right| } =\begin{cases} r^{\left| \lambda\right| }&\text{if } \mathfrak{g}_r=\mathfrak{gl}_r\\ \left( 2r\right) ^{\left|\lambda \right| } &\text{if } \mathfrak{g}_r=\mathfrak{sp}_{2r}

\\ \left( 2r+1\right)^{\left|\lambda \right| } &\text{if } \mathfrak{g}_r=\mathfrak{so}_{2r+1}

\\ \left( 2r\right)^{\left|\lambda \right| } &\text{if } \mathfrak{g}_r=\mathfrak{so}_{2r}

\\\end{cases}.\]  

One also observes a stabilization and independence phenomenon of the tensor multiplicities in large rank. We record it now for later use in the following proposition.

\begin{prop}
\label{si}
Let $\delta,\lambda$ be in $\mathcal{P}_r$. Let $n$ be a positive integer.
\begin{enumerate}
\item If $\delta,n$ are such that $r\geq n\ell\left( \delta\right)$, then one has\[f_{n\,\mu}^{\left( \delta,\mathfrak{g}_r\right) }=f_{n\,\mu}^{\left( \delta,\mathfrak{g}_{r+1}\right)}\qquad \mu \in \mathcal{P}_r.\]

\item If $\delta,\lambda$ are such that $r\geq \ell\left( \lambda\right) +\ell\left( \delta\right)$, then one has\[m_{\lambda\,\mu}^{\left( \delta,\mathfrak{g}_r\right) }=m_{\lambda\,\mu}^{\left( \delta,\mathfrak{g}_{r+1}\right) }\qquad \mu \in \mathcal{P}_r.\]
If $\delta,\lambda,\mu$ are such that $r\geq \ell\left( \lambda\right) +\ell\left( \delta\right)$ and $\left| \mu\right| =\left|\lambda \right| +\left| \delta\right| $, then one has\[m_{\lambda\,\mu}^{\left( \delta,\mathfrak{g}_r\right) }=m_{\lambda\,\mu}^{\left( \delta,\mathfrak{gl}_{r}\right) }.\]
\end{enumerate}

\end{prop}
\subsection{Brief review on the Kashiwara-Nakashima tableaux}
Let $\lambda$ be in $\mathcal{P}_r$ and set $\ell=\ell\left( \lambda\right)$: $\lambda=\left( \lambda_1,\lambda_2,\ldots,\lambda_\ell,0,0,\ldots\right) $. To the irreducible $\mathfrak{g}_r$-module $V\left(\lambda\right)$ is associated its crystal graph $B(\lambda )$, a combinatorial object encoding many informations on $V\left(\lambda\right)$. Here, we only mention the features and properties that are useful for our purposes and refer the reader to the articles~\cite{K},~\cite{KN} and the book~\cite{HK} for a complete exposition.

The graph $B\left(\lambda \right) $ is finite and connected, its number of vertices being equal to the
dimension of $V\left( \lambda \right) $. It is also a colored and oriented graph, the arrows being labelled by the simple roots $S=\left( \alpha_i\right) _{i\in I}$. There exists a weight graduation $\mathrm{wt}:B\left(\lambda\right)\rightarrow P$ such $\mathrm{wt} \left(B\left( \lambda\right)\right)= P\left( \lambda\right)$ and\[v\overset{\alpha_i}{\longrightarrow}{v'}\Longrightarrow \mathrm{wt}\left( v\right)-\mathrm{wt}\left(v'\right)=\alpha_i \qquad v,v'\in B\left( \lambda\right),\,i\in I .\]
There is a unique source vertex $v_0$ and one has $\mathrm{wt}\left(v_0\right)=\lambda$.

When $\mathfrak{g}_r=\mathfrak{gl}_r$, we have $P_+=\mathcal{P}_r$ and the crystal graph $B\left(\lambda\right)$ admits a realization in terms of semistandard Young tableaux of shape $\lambda$ over the alphabet $\mathcal{A}=\left\lbrace 1<\ldots<r\right\rbrace$, that is the fillings of $\lambda$ with letters from the totally ordered set $\mathcal{A}$ such that the rows are weakly increasing as we move to the right and such that the columns are strictly increasing as we go down.

When $\mathfrak{g}_r=\mathfrak{sp}_{2r},\mathfrak{so}_{2r+1}$ or $\mathfrak{so}_{2r}$, there exists a similar parametrization by analogous objects, the so-called Kashiwara-Nakashima $\mathfrak{g}_r$-tableaux. A Kashiwara-Nakashima $\mathfrak{g}_r$-tableau with shape $\lambda$ is a filling of $\lambda $ by letters of the ordered alphabet\[\mathcal{A}=\begin{cases}
\left\lbrace 1<2<\ldots<r<\overline{r}<\ldots<\overline{2}<\overline{1}\right\rbrace&\text{if }\mathfrak{g}_r=\mathfrak{sp}_{2r} \\

\left\lbrace 1<2<\ldots<r<0<\overline{r}<\ldots<\overline{2}<\overline{1}\right\rbrace&\text{if } \mathfrak{g}_r=\mathfrak{so}_{2r+1}\\

\left\lbrace 1<2<\ldots<r-1< r, \overline{r}<\overline{r-1}<\ldots<\overline{2}<\overline{1}\right\rbrace&\text{if } \mathfrak{g}_r=\mathfrak{so}_{2r}\end{cases}\]
according to combinatorial rules depending on $\mathfrak{g}_r$. The set of Kashiwara-Nakashima $\mathfrak{g}_r$-tableaux with shape $\lambda$ contains then all the ordinary semistandard tableaux with shape $\lambda$ along with some generalizations that we do not detail here. One can show that the weight of a Kashiwara-Nakashima $\mathfrak{g}_r$-tableau $T$ of shape $\lambda$ is then the element $\mathrm{wt}\left( T\right)$ of $P$ such that the coordinate on $\varepsilon_i$ is the number of occurrences of the letter $i$ minus the number of occurrences of the letter $\bar{i}$. The source vertex $v_0$ is identified to\[T_0=\begin{array}{cccccc}
1 & \ldots & \ldots & \ldots & 1&\left(\lambda_1 \text{ boxes}\right)  \\ 2 & \ldots & \ldots & 2 &&\left(\lambda_2\text{ boxes}\right)  \\
\vdots& & & &&\\
\ell&\ldots & \ell &&&\left(\lambda_\ell \text{ boxes}\right) \end{array}.\]
Given two adjacent vertices $T\overset{\alpha_i}{\longrightarrow} T'$ in $B\left( \lambda\right) $, they only differ by one letter, the letter in $T$ being replaced in $T'$ by its immediate successor in $\mathcal{A}$. The Weyl character corresponding to $\lambda$ can be computed with these tableaux:\begin{equation}s_\lambda \left( x\right)=\sum_{T\in B\left( \lambda\right) }x^{\mathrm{wt}\left( T\right) }\label{t}.\end{equation}

\section{Markov chains on the finite rank multiplicative graphs}
\label{fr}

In this subsection we fix a nonzero partition $\delta$. We now introduce the multiplicative graph on which the Markov chains we consider are defined.

\subsection{Finite rank multiplicative graphs}
Let $r$ be an integer such that $r\geq \ell \left( \delta\right)$. Define
\[\mathcal{V}^{\left( \delta,\mathfrak{g}_r\right) }=\bigsqcup_{n>0} \mathcal{V}_n^{\left( \delta,\mathfrak{g}_r\right) }\]
where
\[\mathcal{V}_n^{\left( \delta,\mathfrak{g}_r\right) }=\left\lbrace \left( \mu,n\right)  \in \mathcal{P}_r\times \mathbb{N}^* : f_{n\, \mu}^{\left(\delta, \mathfrak{g}_r\right) }>0\right\rbrace\qquad n>0. \]
Given $\left( \mu,n+1\right) $ in $\mathcal{P}_r\times \mathbb{N}^*$, one has\[\left( \mu,n+1\right) \in \mathcal{V}_{n+1}^{\left( \delta,\mathfrak{g}_r\right) } \iff \exists \lambda \in \mathcal{P}_r \quad \left( \lambda,n\right)\in \mathcal{V}_{n}^{\left( \delta,\mathfrak{g}_r\right) }\quad m_{\lambda\, \mu}^{\left(\delta, \mathfrak{g}_r\right) }>0. \]
Let us introduce a set of weighted arrows $\mathcal{E}^{\left( \delta,\mathfrak{g}_r\right) }$ on $\mathcal{V}^{\left(\delta,\mathfrak{g}_r \right) }$: for $\left(\lambda,n\right)$ in $\mathcal{V}_{n}^{\left(\delta, \mathfrak{g}_r\right) }$ and $\left(\mu,n+1\right)$ in $\mathcal{V}_{n+1}^{\left(\delta, \mathfrak{g}_r\right) }$, we set\[\left( \lambda,n\right) \xrightarrow[]{m_{\lambda\, \mu}^{\left(\delta, \mathfrak{g}_r\right) }}\left( \mu,n+1\right)\]if and only if the multiplicity ${m_{\lambda\, \mu}^{\left(\delta, \mathfrak{g}_r\right) }}$ is positive.

\begin{definition}
The graph $\mathcal{G}^{\left(\delta, \mathfrak{g}_r\right) }=\left( \mathcal{V}^{\left(\delta, \mathfrak{g}_r\right) },\mathcal{E}^{\left(\delta, \mathfrak{g}_r\right) }\right) $ is called the multiplicative graph for $\left( \delta,\mathfrak{g}_r\right) $.
\end{definition}

\begin{ex}
\label{examples}
\begin{enumerate}
\item The multiplicative graph for $\left(\tiny{\yng(1)}, \mathfrak{gl}_r\right)  $ is the Young lattice of partitions with at most $r$ parts: $\left(\lambda,n \right) $ in $\mathcal{P}_r\times\mathbb{N}^*$ is a vertex of this graph if and only if $\left| \lambda\right| =n$. The second component of a vertex is then useless. Each arrow has weight $1$. Below are represented  some levels of the multiplicative graph for $\left(\tiny{\yng(1)}, \mathfrak{gl}_{2}\right)  $:
\[
{\tiny
	\begin{array}{ccccccccccccccccc}

	&&&&&&&& \yng( 1) &&&&&&&&\\
	
	&&&&&&& \swarrow && \searrow &&&&&&&\\
	
	&&&&&& \yng( 2) &&&& \yng( 1,1) &&&&&&\\
	
	&&&&&\swarrow &&\searrow&&\swarrow&&&&&&&\\
	
	&&&&\yng( 3) &&&&\yng(2,1) &&&&&&& \\
	
	&&&\swarrow &&\searrow&&\swarrow&\downarrow&&&&&&&&\\
	
	&&\yng( 4)&&&&\yng(3,1)&&\yng(2,2)&&&&&&&&
	
	\end{array}}
\]
\item The multiplicative graph for $\left(\tiny{\yng(1)},\mathfrak{sp}_{2r} \right)$ is the Pascalized graph of the previous one: starting from the multiplicative graph for $\left(\tiny{\yng(1)}, \mathfrak{gl}_r\right)$, we complete the level $n+1$ by reflecting the level $n-1$ with respect to the level $n$, with the corresponding arrows. Each arrow has weight $1$. Then $\left(\lambda,n \right) $ in $\mathcal{P}_r\times\mathbb{N}$ is a vertex of this graph if and only if $\left|\lambda\right|\leq n $ and $\left|\lambda\right|= n \mod 2$. Below are represented  some levels of the multiplicative graph for $\left(\tiny{\yng(1)},\mathfrak{sp}_{4}\right)$:
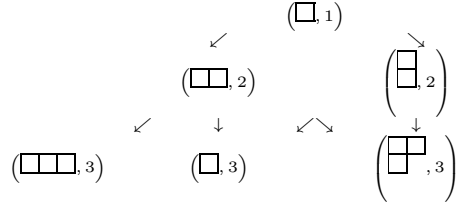
\[
{\tiny
	\begin{array}{ccccccccccccccccc}

&&&&&&&& \left( \yng( 1),1\right)  &&&&&&&&\\

&&&&&&& \swarrow && \searrow &&&&&&&\\

&&&&&&& \left( \yng( 2),2\right)  && \left( \yng( 1,1),2\right)  &&&&&&&\\

&&&&&&\swarrow&\downarrow&\swarrow\searrow &\downarrow&&&&&&&\\

&&&&  &\left(\yng( 3),3\right)&&\left( \yng(1),3\right) & &\left( \yng(2,1),3\right) &&&&&& \\
\end{array}}\]
\item Let $\left(\lambda,n \right)$ be in $\mathcal{P}_r\times\mathbb{N}^*$. If it is a vertex of the multiplicative graph for $\left( \tiny{\yng(2)},\mathfrak{so}_{2r+1}\right)$, then it must satisfy $\left|\lambda \right|\leq \left| \tiny{\yng(2)}\right|\times n=2n$ and $\left|\lambda \right|=0 \mod 2$. This condition is not sufficient since $\left( \tiny{\yng(2,1,1)},2\right)$ is not a vertex of the multiplicative graph for $\left(\tiny{\yng(2)},\mathfrak{so}_{7} \right)$. Beside this, there are adjacent vertices with the same number of boxes and some weights are greater than $1$: in the multiplicative graph for $\left(\tiny{\yng(2)},\mathfrak{so}_{7} \right) $, one has the weighted arrow \[\left( \tiny{\yng( 3,1)},2\right)  \xrightarrow[]{2}\left( \tiny{\yng( 3,1)},3\right)  .\] 
In particular this last multiplicative graph is not the Pascalized version of a braching graph of type $A$.
\end{enumerate}
\end{ex}

\subsection{Markov chains on the finite rank multiplicative graphs}
Let $r$ be an integer such that $r\geq \ell\left( \delta\right) $. In this paragraph, we introduce the normalized Weyl characters and then use them to define a stochastic matrix on the adjacent vertices of $\mathcal{G}^{\left( \delta,\mathfrak{g}_r\right) }$.

Let $\lambda$ be in $\mathcal{P}_r$. We have\[x^{-\lambda}s_\lambda^{\mathfrak{g}_r} \left({x}\right) =\sum_{\omega\in P^{\mathfrak{g}_r}\left( \lambda\right)  }K_{\lambda\,\omega}^{\mathfrak{g}_r} {x}^{-\left(\lambda-\omega \right) }.\]
If $\omega$ is in $P^{\mathfrak{g}_r}\left( \lambda\right) $, then there exists $T$ in $B^{\mathfrak{g}_r}\left( \lambda\right) $ such that $\mathrm{wt}\left( T\right) =\omega$ and there is a path in $B^{\mathfrak{g}_r}\left( \lambda\right) $ starting at $T_0$ and ending at $T$. Then $\lambda-\omega=\mathrm{wt}\left( T_0\right) -\mathrm{wt}\left( T\right) $ must be a linear combination of the simple roots $S=\left( \alpha_i\right) _{i\in I}$ with nonnegative integer coefficients: more precisely, the coefficient of $\alpha_i$ in $\mathrm{wt}\left( T_0\right) -\mathrm{wt}\left( T\right)$  is equal to the number of browsed arrows $\xrightarrow[]{\alpha_i}$ in any path from $T_0$ to $T$. Using the family $\Omega^{\vee}=\left( \omega_i^{\vee}\right) _{i\in I}$, one can also write\[\lambda-\omega=\sum_{i\in I}\left<\lambda-\omega,\omega_i^{\vee} \right>\alpha_i.\] Define a family of variables $\left( y_i\right) _{i\in I}$ by setting $y_i=x^{-\alpha_i}$. If we introduce ${y}^{\left[\omega\right]}=\prod_{i\in I} y_i^{\left<\omega,\omega_i^{\vee} \right> }$ for $\omega$ in $L$, then we can write\[\sum_{\omega\in P^{\mathfrak{g}_r}\left( \lambda\right)  }K_{\lambda\,\omega}^{\mathfrak{g}_r} {x}^{-\left(\lambda-\omega \right) }=\sum_{\omega\in P^{\mathfrak{g}_r}\left( \lambda\right) }K_{\lambda \, \omega}^{\mathfrak{g}_r}{y}^{\left[\lambda-\omega \right]}.\]
Thus $S_\lambda ^{\mathfrak{g}_r}\left({y}\right)=\sum_{\omega\in P^{\mathfrak{g}_r}\left( \lambda\right) }K_{\lambda \, \omega}^{\mathfrak{g}_r}{y}^{\left[\lambda-\omega \right]}$ is a polynomial in $\left( y_i\right) _{i\in I}$ with nonnegative integer coefficients. Given a finite sequence of positive reals ${\theta}=\left(\theta_i\right)_{i\in I}$ and $\omega$ in $L$, the real number obtained by letting $y_i=\theta_i$ in ${y}^{\left[ \omega\right] }$ (respectively $S_\lambda^{\mathfrak{g}_r}\left({y}\right)$), is denoted by $\theta^{\left[ \omega\right] }$ (respectively $S_\lambda^{\mathfrak{g}_r}\left(\theta\right)$). When $b$ is a positive real, we write $S_\lambda^{\mathfrak{g}_r} \left({ b,b,\ldots,b }\right)$ for the real number obtained by replacing each $y_i$ by $b$ in $S_\lambda^{\mathfrak{g}_r} \left({y}\right)$. In other words: \begin{equation}S_\lambda^{\mathfrak{g}_r} \left({b,b,\ldots,b}\right) =\sum_{\omega\in P^{\mathfrak{g}_r}\left( \lambda\right) } K_{\lambda\,\omega}^{\mathfrak{g}_r}b^{\left\langle \lambda-\omega,\rho^{\vee} \right\rangle  }\end{equation}
where $\rho^{\vee}=\sum_{i\in I }\omega_i^{\vee}$.

For any finite sequence of positive reals ${\theta}=\left(\theta_i\right)_{i\in I}$, we define a matrix $\Pi_\theta^{\left(\delta,\mathfrak{g}_r \right) }$ on $\mathcal{G}^{\left( \delta,\mathfrak{g}_r\right) }$:  for any $\left( \lambda,n\right) ,\left(\mu ,n+1\right) $ in $\mathcal{P}_r \times \mathbb{N}^*$ such that
\[\left( \lambda,n\right) \xrightarrow{m_{\lambda\,\mu}^{\left(\delta,\mathfrak{g}_r\right) }}\left( \mu,n+1\right),\]set \[\Pi_\theta^{\left(\delta,\mathfrak{g}_r \right) }\left(\left( \lambda,n\right),\left( \mu,n+1\right) \right)=m_{\lambda\,\mu}^{\left( \delta,\mathfrak{g}_r\right) }\frac{S_{\mu}^{\mathfrak{g}_r} \left(\theta\right) }{S_{\lambda}^{\mathfrak{g}_r}\left(\theta \right)S_{\delta}^{\mathfrak{g}_r}\left(\theta \right)}{ \theta^{\left[\lambda+\delta-\mu\right]}}.\]
It is a stochastic matrix because~\eqref{pie} implies\[s_\lambda^{\mathfrak{g}_r} \left( {x}\right)s_\delta^{\mathfrak{g}_r} \left( {x}\right)=\sum_{\mu \in \mathcal{P}_r } m_{\lambda\,\mu}^{\left( \delta,\mathfrak{g}_r\right) }s_\mu ^{\mathfrak{g}_r}\left({ x}\right)\qquad \lambda \in \mathcal{P}_r.\]For every fixed parameter $\theta$, we thus obtain a Markov chain $\mathcal{M}_\theta^{\left( \delta,\mathfrak{g}_r\right)}$ associated to the pair $\left( \mathcal{G}^{\left( \delta,\mathfrak{g}_r\right) },\Pi_\theta^{\left(\delta,\mathfrak{g}_r \right) }\right) $.

\section{Markov chains on the infinite rank multiplicative graphs}
\label{irmc}

Let us fix once for all a nonzero partition $\delta$. Let $X=\left(\mathfrak{g}_r\right)_r$ be $A,C,B$ or $D$.

\subsection{Infinite rank multiplicative graphs}
\label{irg}

For fixed $\left(\mu,n \right) $ in $\mathcal{P}_\infty \times \mathbb{N}^*$, the first point of Proposition~\ref{si} shows that the sequence of integers $\left(f_{n\,\mu}^{\left( \delta,\mathfrak{g}_r\right) } \right)_r$ eventually becomes constant. Denote by $f_{n\,\mu}^{\left( \delta,X\right) }$ this stabilized multiplicity. So $f_{n\,\mu}^{\left( \delta,X\right) }$ is positive exactly means that, for all $r$ sufficiently large, $\left(\mu ,n\right)$ is an element of $\mathcal{V}_n^{\left(\delta,\mathfrak{g}_r\right) }$. This naturally leads us to define a limit set of vertices by setting \[\mathcal{V}^{\left(\delta,X\right) }=\bigsqcup_{n>0} \mathcal{V}_n^{\left( \delta,X\right) }\]
where
\[\mathcal{V}_n^{\left( \delta,X\right) }=\left\lbrace \left( \mu,n\right)  \in \mathcal{P}_\infty\times \mathbb{N}^* : f_{n\, \mu}^{\left(\delta, X\right) }>0\right\rbrace\qquad n>0. \]
For fixed $\lambda,\mu$ in $\mathcal{P}_\infty$, the second point of Proposition~\ref{si} shows that the sequence of integers $\left(m_{\lambda\, \mu}^{\left( \delta,\mathfrak{g}_r \right)} \right) _r$ eventually becomes constant. Denote by $m_{\lambda\, \mu}^{\left( \delta,X\right) }$ this stabilized multiplicity. Let us introduce a set of weighted arrows $\mathcal{E}^{\left(\delta, X\right) }$ on $\mathcal{V}^{\left(\delta, X\right) }$: for $\left( \lambda,n\right)$ in $\mathcal{V}_{n}^{\left(\delta,X\right) }$ and $\left( \mu,n+1\right) $ in $\mathcal{V}_{n+1}^{\left(\delta, X\right) }$, we set
\[\left( \lambda,n\right) \xrightarrow[]{m_{\lambda\, \mu}^{\left(\delta, X\right) }}\left( \mu,n+1\right) \]
if and only if the stabilized multiplicity $m_{\lambda\, \mu}^{\left( \delta,X\right) }$ is positive.

\begin{definition}
	The graph $\mathcal{G}^{\left(\delta, X\right) }=\left( \mathcal{V}^{\left(\delta, X\right) },\mathcal{E}^{\left(\delta, X\right) }\right) $ is called the multiplicative graph for $\left( \delta,X\right) $.
\end{definition}

\begin{ex}
\begin{enumerate}
\item The multiplicative graph for $\left(\delta,A \right) $ is the classical Young lattice. Each arrow has weight $1$. Here is represented some levels: 

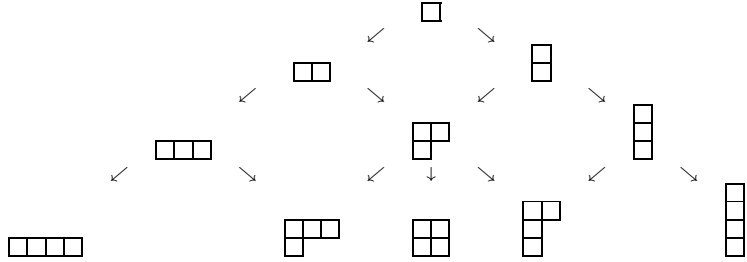
\[
{\tiny
	\begin{array}{ccccccccccccccccc}

	&&&&&&&& \yng( 1) &&&&&&&&\\
	
	&&&&&&& \swarrow && \searrow &&&&&&&\\
	
	&&&&&& \yng( 2) &&&& \yng( 1,1) &&&&&&\\
	
	&&&&&\swarrow &&\searrow&&\swarrow&&\searrow&&&&&\\
	
	&&&&\yng( 3) &&&&\yng(2,1) &&&&\yng( 1,1,1)&&& \\
	
	&&&\swarrow &&\searrow&&\swarrow&\downarrow&\searrow&&\swarrow&&\searrow&&&\\
	
	&&\yng( 4)&&&&\yng(3,1)&&\yng(2,2)&&\yng(2,1,1)&&&&\yng(1,1,1,1)&&
	
	\end{array}}
\]

\item The multiplicative graph for $\left( \tiny{\yng(1)},C\right) $ is the Pascalized of the previous one.
\item The multiplicative graph for $\left( \tiny{\yng(2)},B\right) $ is not a Pascalized version of a multiplicative graph of type $A$.
\end{enumerate}
\end{ex}

\subsection{A common stochastic matrix for the infinite rank multiplicative graphs} 
\label{irm}

Let $\theta=\left( \theta_i\right) _{i\geq 1}$ be a sequence of positive reals and set \[\theta_{\left[ r\right] }=\left(\theta_i \right)_{i\in I}=\begin{cases} \left( \theta_1,\theta_2,\ldots,\theta_{r-1}\right) &\text{if } X=A\\
\left( \theta_1,\theta_2,\ldots,\theta_{r}\right)&\text{if } X=C,B,D
\end{cases}\] for any positive integer $r$.

Our aim is now, under a suitable hypothesis on $\theta$, to establish
the pointwise convergence of the sequence of matrices $\left(\Pi_{\theta_{\left[ r\right] }}^{\left( \delta,\mathfrak{g}_r\right) } \right)_r$ on the set $\mathcal{P}_\infty\times \mathbb{N}^*$ and also identify the limit.
 
We start by showing that the sequence $\left(S_\lambda^{\mathfrak{g}_r} \left( \theta_{\left[ r\right] } \right) \right)_r $ is convergent for any $\lambda$ in $\mathcal{P}_\infty$. Our method consists in showing this result for $X=A$ and then use it in the other cases.

\subsubsection{Case $X=A$}
\label{A}
Let $\lambda$ be in $\mathcal{P}_\infty$ and set $\ell=\ell\left( \lambda\right) $: $\lambda=\left( \lambda_1,\lambda_2,\ldots,\lambda_\ell,0,0,\ldots\right)$. Assume $X=A$. Our first ingredient is the following proposition. Our proof is directly inspired by the proof of the dimension formula given in~\cite{FH}.

\begin{prop}\label{pr}
Let $b$ be a positive real. One has
\[S_\lambda ^{\mathfrak{gl}_r}\left(b,b,\ldots,b\right)=\prod_{1\leq i<j\leq \ell}\frac{1-b^{\lambda_i -\lambda_j +j-i}}{1-b^{j-i}}\prod_{1\leq i\leq \ell}\, \prod_{\ell<j\leq r}  \frac{1-b^{\lambda_i+j-i}}{1-b^{j-i}}\qquad r>\ell.\]
\end{prop}
\textit{Proof. }In this proof, we omit to include the parameter $\mathfrak{gl}_r$ in our notation. Let $r$ be an integer such that $r>\ell$. Let $z$ be a variable. For any $\omega_0$ in $L=\bigoplus_{i=1}^r\mathbb{Z}\varepsilon_i$, introduce now the morphism $\phi_{\omega_0}:\mathbb{C}\left[L\right] \rightarrow   \mathbb{C}\left[\left[z\right] \right]$ defined by $\phi_{\omega_0}\left( {x}^{\omega}\right)=\exp\left({\left< \omega,\omega_0\right> z}\right)$ for any $\omega$ in $L$.

Firstly, we have\[S_\lambda \left(b,b,\ldots,b\right)=\sum_{\omega\in P\left( \lambda\right)  }K_{\lambda\,\omega}\exp \left(\left<\lambda-\omega,\rho\right>\ln b \right)\]since $\rho^{\vee}=\rho$ when $X=A$. Observe that $S_\lambda \left(b,b,\ldots,b\right)$ is obtained by letting $z=-\ln b$ in $\phi_{\rho}\left({x}^{-\lambda}s_\lambda \left({x}\right) \right)$. 

Secondly, the Weyl character formula gives \[{x}^{-\lambda}s_\lambda\left({x}\right)={x}^{-\lambda}\frac{a_{\lambda+\rho}\left({x}\right) }{a_{\rho}\left( {x}\right)}.\]
 For any $\omega_0,\omega$ in $L$, observe that we have\[\begin{aligned}\phi_{\omega_0} \left(a_{\omega}\left({x}\right) \right) =\sum_{\sigma \in W}\epsilon\left( \sigma\right) \phi_{\omega_0} \left(x^{\sigma \omega}\right) &=\sum_{\sigma \in W}\epsilon\left( \sigma\right) \exp\left(\left< \sigma \omega,\omega_0\right>z  \right)\\&=\sum_{\sigma \in W}\epsilon\left( \sigma^{-1}\right) \exp\left(\left< \omega,\sigma^{-1}\omega_0\right>z  \right)\\&= \phi_{\omega}\left( a_{\omega_0}\left({x}\right)\right).\end{aligned}\] Then \[\begin{aligned}\phi_{\rho}\left(x^{-\lambda}a_{\lambda+\rho}\left({x}\right) \right)&= \phi_{\rho}\left(x^{-\lambda}\right)  \phi_{\lambda+\rho}\left(a_{\rho}\left( x\right)  \right)\\&=\exp\left(\left\langle -\lambda,\rho\right\rangle z \right)\phi_{\lambda+\rho}\left(x^{\rho}\prod_{\alpha\in R_+}\left(1-x^{-\alpha} \right) \right)\\
&=\exp\left( \left<\rho,\rho\right>z\right)\prod_{\alpha \in R_+}\left( 1-\exp\left(-\left< \alpha,\lambda+\rho\right> z \right) \right) \end{aligned}\]and
\[\phi_{\rho}\left(a_{\rho}\left({x}\right) \right)=\exp\left( \left<\rho,\rho\right>z\right)\prod_{\alpha \in R_+}\left( 1-\exp\left(-\left<\alpha,\rho\right> z \right) \right) .\]
This shows\[S_\lambda \left(b,b,\ldots,b\right)=\prod_{\alpha\in R_+} \frac{1-b^{\left\langle \lambda+\rho,\alpha\right\rangle }}{1-b^{\left\langle \rho,\alpha\right\rangle }}.\]
As $X=A$, we have\[R_+=\left\lbrace \varepsilon_i-\varepsilon_j : 1\leq i<j\leq r\right\rbrace \]and $\rho=\sum_{i=1}^r \left(r-i\right) \varepsilon_i$. Reminding that $\lambda_j=0$ for any integer $j$ such that $j>\ell$, the result follows.$\Box$\\

\begin{lem}
\label{p} For any real $t$ in $\left(0,1\right) $ and for any positive integers $a$, the infinite product $\prod_{n>0}\frac{1-t^{a+n}}{1-t^{n}}$ is convergent.
\end{lem}
\textit{Proof.} We have indeed $\ln \left( \frac{1-t^{a+n}}{1-t^n}\right)\underset{n\to\infty}{\sim} \left( 1-t^a\right) t^n$.$\Box$\\

\begin{theo}
\label{n}
 Assume $\theta$ is a bounded sequence of positive reals such that $\sup \theta <1$. The sequence $\left(S_{\lambda}^{\mathfrak{gl}_r}\left(\theta_{\left[ r\right] }\right)\right) _{r}$ is convergent. The limit is denoted by $S_{\lambda}^{A}\left(\theta\right)$ in the sequel.
\end{theo}
\textit{Proof.} Any semistandard tableau on the alphabet $\left\lbrace 1<2<\ldots<r\right\rbrace $ is also a semistandard tableau on the alphabet $\left\lbrace 1<2<\ldots<r<r+1\right\rbrace $, with same weight. Then Identity~\eqref{t} implies that the sequence $\left(S_\lambda^{\mathfrak{gl}_r} \left(\theta_{\left[ r\right] } \right) \right) _r$ is weakly increasing.

Recall that, for any positive integer $r$ such that $r\geq \ell$, the polynomial $S_\lambda^{\mathfrak{gl}_r}\left(y\right)$ has nonnegative integer coefficients. We have thus\[S_\lambda^{\mathfrak{gl}_r}\left( \theta_{\left[ r\right] }\right)\leq S_\lambda^{\mathfrak{gl}_r} \left(b,b\ldots,b\right)\qquad r\geq \ell\]where $b=\sup \theta$. Since $b<1$, Proposition~\ref{pr} and Lemma~\ref{p} show that the sequence $\left(S_\lambda^{\mathfrak{gl}_r} \left(b,b\ldots,b\right) \right)_r $ is bounded and so is $\left(S_{\lambda}^{\mathfrak{gl}_r}\left(\theta_{\left[ r\right] }\right)\right) _{r}$. It is thus convergent.$\Box$\\

\subsubsection{Case $X=C$, $B$, or $D$}
\label{CBD}
Let $\lambda$ be in $\mathcal{P}_\infty$ and set $\ell=\ell\left( \lambda\right) $. Assume $X=C,B$, or $D$. Assume $\theta$ is bounded and $\sup \theta <1$. Set $b=\sup \theta$.
\begin{prop}  We have
\[0\leq S_\lambda^{\mathfrak{g}_r}\left( \theta_{\left[ r\right]} \right)-S_\lambda ^{\mathfrak{gl}_r} \left( \theta_{\left[ r\right]}\right)\leq \begin{cases} b^{r-\ell+1} \left( 2r\right) ^{\left| \lambda\right|}&\text{if }\mathfrak{g}_r=\mathfrak{sp}_{2r}\\
b^{r-\ell} \left( 2r+1\right) ^{\left| \lambda\right|}&\text{if }\mathfrak{g}_r=\mathfrak{so}_{2r+1}\\
b^{r-\ell+2} \left( 2r\right) ^{\left| \lambda\right|}&\text{if }\mathfrak{g}_r=\mathfrak{so}_{2r}
\end{cases}\qquad r\geq \ell.\]
\end{prop}
\textit{Proof.} Suppose $X=C$, the arguments are similar in the other cases. Let $r$ be such that $r\geq\ell$. We have\[S_\lambda ^{\mathfrak{g}_r}\left(\theta_{\left[ r\right]}\right)=\sum_{T\in B^{\mathfrak{g}_r}\left(\lambda\right) }\theta_{\left[ r\right]}^{\left[ \lambda-\mathrm{wt}\left(T\right) \right] }.\] Since $B^{\mathfrak{g}_r}\left( \lambda\right)$ contains $B^{\mathfrak{gl}_r}\left(\lambda\right)$ and since the simple roots $\alpha_1,\alpha_2,\ldots,\alpha_{r-1}$ are the same for $\mathfrak{g}_r$ and $\mathfrak{gl}_r$, the left inequality is clear.

Now let $T$ be in $B^{\mathfrak{g}_r}\left( \lambda\right)\setminus B^{\mathfrak{gl}_r}\left( \lambda\right)$: there is at least one letter of $T$ which is in $\left\lbrace \overline{r},\overline{r-1},\ldots,\overline{1}\right\rbrace$. A path in $B^{\mathfrak{g}_r}\left( \lambda\right)$ starting at $T_0$ and ending at $T$ must browse arrows $\overset{\alpha_\ell}{\longrightarrow},\ldots,\overset{\alpha_r}{\longrightarrow}$, each of them at least one time, in order to change one letter $\ell$ into an element of $\left\lbrace \overline{r},\overline{r-1},\ldots,\overline{1}\right\rbrace$. So, looking at the weight graduation, we have to substract at least one time each of the simple roots $\alpha_\ell,\ldots,\alpha_r$ in order to get $\mathrm{wt}\left( T\right) $ from $\mathrm{wt}\left( T_0\right) $. Since $b<1$, one can write then\[\begin{aligned}S_\lambda ^{\mathfrak{g}_r}\left(\theta_{\left[ r\right] }\right)-S_\lambda ^{\mathfrak{gl}_r}\left(\theta_{\left[ r\right]}\right)&=\left(\prod_{i=\ell}^r \theta_i\right)\sum_{T\in B^{\mathfrak{g}_r}\left( \lambda\right)\setminus B^{\mathfrak{gl}_r}\left( \lambda\right) }\theta_{\left[ r\right] }^{\left[ \lambda-\mathrm{wt}\left(T\right)-\sum_{i=\ell}^r\alpha_i\right] }\\ &\leq  b^{r-\ell+1} \sum_{T\in B^{\mathfrak{g}_r}\left( \lambda\right)\setminus B^{\mathfrak{gl}_r}\left( \lambda\right) }1\\&\leq b^{r-\ell+1} \dim V^{\mathfrak{g}_r} \left( \lambda\right)\\&\leq  b^{r-\ell+1} \left( 2r\right) ^{\left| \lambda\right| }.\end{aligned}\]
$\Box$\\

\begin{cor}
\label{nA}
We have \[S_\lambda^{\mathfrak{g}_r} \left( \theta_{\left[ r\right] }\right)\xrightarrow[r\to\infty]{}S_\lambda^{A}\left( \theta\right).\]
\end{cor}
\textit{Proof.} As $b=\sup \theta <1$, we have\[b^{r} r^{\left| \lambda\right|}\xrightarrow[r\to \infty]{}0\] and the result is then a straightforward consequence of the previous proposition.$\Box$\\

\subsubsection{Main result}
\label{main}
Assume $X=C,B$, or $D$. Assume $\theta$ is bounded and $\sup \theta<1$.
\begin{lem}
\label{l}
Let $\lambda,\mu$ be in $\mathcal{P}_\infty$. If $\left|\mu\right|=\left| \lambda\right| +\left| \delta\right|$, then the sequence $\left( \theta_{\left[ r\right] }^{\left[ \lambda+\delta-\mu\right] }\right) _r$ eventually becomes constant. If $\left|\mu\right| <\left| \lambda\right| +\left| \delta\right|$, then \[{\theta_{\left[ r\right] }^{\left[ \lambda+\delta-\mu\right] } }\xrightarrow[r\to \infty]{}0.\]
\end{lem}
\textit{Proof.} Suppose $X=B$, the arguments are similar in the other cases. Let us fix a positive integer $\ell$ such that $ \lambda+\delta-\mu$ is in $\bigoplus_{i=1}^{\ell} \mathbb{Z}\varepsilon_i$. As the family $\left(\varepsilon_1-\varepsilon_{2},\ldots,\varepsilon_{\ell-1}-\varepsilon_{\ell},\varepsilon_{\ell} \right) $ is a basis of $\bigoplus_{i=1}^{\ell} \mathbb{Z}\varepsilon_i$, there exists integers $a_1,\ldots,a_{\ell}$ such that \[{\lambda+\delta }-\mu=\sum_{i=1}^{\ell-1} a_i \left( \varepsilon_i-\varepsilon_{i+1}\right) +a_{\ell} \varepsilon_{\ell}.\] Now observe that $a_{\ell}=\left| \lambda\right| +\left| \delta\right|-\left|\mu\right|$. Recall that $\varepsilon_{i}-\varepsilon_{i+1}=\alpha_i$ for any integer $i$ such that $1\leq i <\ell$ and $\varepsilon_{\ell}=\sum_{i=\ell}^{r} \alpha_i$ for any integer $r$ such that $r\geq \ell$. We deduce\[\theta_{\left[ r\right] }^{\left[\lambda+\delta -\mu \right]}= \prod_{i=1}^{\ell-1} \theta_i^{a_i} \left( \prod_{i=\ell}^{r} \theta_i\right)^{a_\ell}\qquad r\geq \ell.\]
If $a_{\ell}=\left| \lambda\right| +\left| \delta\right|-\left|\mu\right|=0$, then\[\theta_{\left[ r\right] }^{\left[\lambda+\delta -\mu \right]}= \prod_{i=1}^{\ell-1} \theta_i^{a_i}\qquad r\geq\ell.\]
If $a_{\ell}=\left| \lambda\right| +\left| \delta\right|-\left|\mu\right|>0$, then the assumption $\sup \theta<1$ implies\[\theta_{\left[ r\right] }^{\left[\lambda+\delta -\mu \right]}\leq \prod_{i=\ell}^{r} \theta_i\qquad r\geq \ell\]
and
\[\prod_{i=\ell}^{r} \theta_i\xrightarrow[r\to\infty]{} 0.\]$\Box$\\

Now, we state our main result for which $X$ can be $A,C,B$ or $D$.
\begin{theo}  
\label{mainTh}
For any $\left( \lambda,n\right) ,\left(\mu ,n+1\right) $ in $\mathcal{P}_\infty \times \mathbb{N}^*$ such that
\[\left( \lambda,n\right) \xrightarrow{m_{\lambda\,\mu}^{\left(\delta,X\right) }}\left( \mu,n+1\right),\] one has
\[\Pi_{\theta_{\left[ r\right] }}^{\left(\delta,\mathfrak{g}_r \right) } \left( \left( \lambda,n\right) ,\left( \mu,n+1\right)\right) \xrightarrow[r\to\infty]{}\begin{cases} m_{\lambda\,\mu}^{\left( \delta,A\right) }\frac{S_{\mu}^{A} \left(\theta\right) }{S_{\lambda}^{A}\left(\theta \right)S_{\delta}^{A}\left(\theta \right)}{ \theta^{\left[\lambda+\delta-\mu\right] } }&\text{if } \left| \mu\right| =\left|\lambda \right| +\left|\delta \right| \\ 
0&\text{if} \left| \mu\right|<\left|\lambda \right| +\left|\delta \right|\end{cases}\]
\end{theo}
\textit{Proof.} Proposition~\ref{si} and Corollary~\ref{nA} show the convergence of the sequences\[\left(m_{\lambda\,\mu}^{\left(\delta,\mathfrak{g}_r \right)} \right), \left(S_{\delta}^{\mathfrak{g}_r }\left(\theta_{\left[ r\right] } \right) \right)_r,\left(S_{\lambda}^{\mathfrak{g}_r }\left(\theta_{\left[ r\right] } \right) \right)_r,\left(S_{\mu}^{\mathfrak{g}_r }\left(\theta_{\left[ r\right] } \right) \right)_r\]to the positive reals\[m_{\lambda\,\mu}^{\left( \delta,X\right)},S_{\delta}^{A}\left(\theta \right),S_{\lambda}^{A}\left(\theta \right),S_{\mu}^{A}\left(\theta \right),\] respectively. When $\left| \mu\right| =\left|\lambda \right| +\left|\delta \right|$, Proposition~\ref{si} also mentions that we have $m_{\lambda\,\mu}^{\left( \delta,X\right) }=m_{\lambda\,\mu}^{\left( \delta,A\right) }$. Then Lemma~\ref{l} gives the expected result.$\Box$\\

\section{Further results}
\label{SecFR}

\subsection{Connection with the generalized Pitman transform}

It was shown in \cite{LLP2} that the Markov chains defined in \ref{irmc}
coincide with some random Littelmann paths conditioned to stay in Weyl
chambers and can also be obtained from non conditioned such random paths by
applying the generalized Pitman transform introduced in \cite{BBOC1}.\ For
$\mathfrak{gl}_{r}$, this generalized Pitman transform can also be described
from the Schensted insertion algorithm on semistandard tableaux.\ Similar
algorithms on Kashiwara-Nakashima tableaux also exist for the other classical
types in finite rank based on Kashiwara crystal basis theory.

The generalized Pitman transform essentially associates to each concatenation
of Littelmann paths for the representation $V(\delta)$ considered its
corresponding highest weight path which entirely lies in the Weyl chamber.
This construction cannot be directly generalized in infinite rank because the
relevant representations to consider are not then of highest weight.
Nevertheless the crystal basis construction can be adapted together with the
insertion algorithms on tableaux (see \cite{L}). It is interesting to observe
that the same kind of stabilization phenomenons then appear in large rank:
the generalized Pitman transforms admit a limit in large rank and the four
limits so obtained in each classical types are essentially the same.

\subsection{Extremal harmonic functions on the limit multiplicative graphs}

The extremal harmonic functions on the multiplicative graphs $\mathcal{G}%
^{(\delta,\mathfrak{g}_{r})}$ with $\mathfrak{g}_{r}$ of type $X$ have
been characterized in \cite{LT} in terms of the functions $S_{\lambda}%
(\theta_{\lbrack r]})$ introduced in \S \ 3.2.\ It follows from our main Theorem \ref{mainTh}
that the extremal harmonic functions on the limit multiplicative graphs
$\mathcal{G}^{(\delta,X)}$ only depend on the partition $\delta$ considered
and not on the type $X$ considered. It thus suffices to consider the extremal
harmonic functions on $\mathcal{G}^{(\delta,A)}$ which can be expressed in
terms of Schur functions by refining the work of Kerov and Vershik on the
Young lattice (see \cite{LT}).

\end{document}